\documentclass[12pt]{article}
\usepackage{amsfonts}
\makeatletter
 
\@addtoreset{equation}{section}
\def\theequation{\thesection.\arabic{equation}}
 
\makeatother
\renewcommand\theequation{\thesection.\arabic{equation}}

\newtheorem{thm}{Theorem}[section]
\newtheorem{lemma}[thm]{Lemma}

\begin{document}
\def\E{{\mathbb E}}
\def\P{{\mathbb P}}
\def\R{{\mathbb R}}
\def\Z{{\mathbb Z}}
\def\V{{\mathbb V}}
\def\N{{\mathbb N}}
\def\X{{\cal X}}
\def\G{{\cal G}}
\def\n{{\bf n}}
\def\m{{\bf m}}
\def\sqr{\vcenter{
         \hrule height.1mm
         \hbox{\vrule width.1mm height2.2mm\kern2.18mm\vrule width.1mm}
         \hrule height.1mm}}                  
\def\square{\ifmmode\sqr\else{$\sqr$}\fi}
\def\one{{\bf 1}\hskip-.5mm}
\def\liml{\lim_{L\to\infty}}
\def\given{\ \vert \ }
\def\ze{{\zeta}}
\def\la{{\lambda}}
\def\ga{{\gamma}}
\def\reff#1{(\ref{#1})}
\def\w{{\omega}}
\def\a{\alpha}
\def\T{\Theta}
\def\vep{\varepsilon}
\def\LL{{\cal L}}
\def\L{\Lambda} 
\def\S{{\cal S}}
\def\t{\theta}
\def\comp{\leftrightarrow}
\def\W{W}
\def\ml{{\mu_\L}}
\def\bx{{\bf X}}
\def\e{{\underline\eta}}
\def\M{M}
\def\r{r}
\def\rr{R}
\def\liml{\lim_{\L\to\LL}}
\def\lime{\lim_{\vep\to0}}
\def\A{{\bf A}}
\def\B{{\bf B}}
\def\C{{\bf C}}
\def\D{{\bf D}}
\def\MM{{\bf m}}
\def\proof{\noindent{{\bf Proof. }}}

\title{Blocking measures for asymmetric exclusion processes\\
via coupling } 
\date{}

\vspace{8mm}
\author{
\vspace{2mm}
P. A. Ferrari\thanks {\tt email: \tt pablo@ime.usp.br}\\
{\it Universidade de S\~{a}o Paulo} \\
\vspace{5mm}
{\it IME USP, Caixa Postal 66281, 05315-970
 - S\~{a}o Paulo, BRAZIL}\\
\vspace{2mm}
J. L. Lebowitz\thanks{Also Department of Physics, Rutgers.  email:
{\tt lebowitz@math.rutgers.edu}}
 and E. Speer\thanks{email: {\tt speer@math.rutgers.edu}} \\
 {\it Department of Mathematics, Rutgers University}\\
 {\it New Brunswick, NJ 08903 USA.}}

\maketitle

 \abstract
 We give sufficient conditions on the rates of two asymmetric exclusion
processes such that the existence of a blocking invariant measure for the
first implies the existence of such a measure for the second.  The main tool
is a coupling between the two processes under which the first dominates the
second in an appropriate sense.  In an appendix we construct
a class of processes for which the existence of a blocking measure can be
proven directly; these are candidates for comparison processes in
applications of the main result.

 \vskip0.5truein

 \noindent
 {\bf Key words:}
 {\bf AMS Classification:} Primary:
  60K35 60D05, Secondary: 60G55

\endabstract

\newpage

\setlength{\baselineskip}{21pt}

\section{Introduction}
\label{Introduction}

We consider the exclusion process $\eta_t$ on $\{0,1\}^\Z$ with generator $L$
given by
\begin{equation}
  \label{1}
  Lf(\eta) = \sum_{x\in\Z}\sum_{y\in\Z} p(x,y;\eta)[f(\eta^{x,y}) - f(\eta)]\,.
 \end{equation}
 Here $f$ is a continuous function on $\{0,1\}^\Z$ (with the product
topology) and for $x,y,z\in\Z$, 
 \begin{equation}
\eta^{x,y}(z)=\cases{\eta(x),& if $z=y$,\cr \eta(y),& if $z=x$,\cr
      \eta(z),& otherwise.\cr}
 \end{equation}
 The jump rate of particles from $x$ to $y$ in configuration $\eta$,
$p(x,y;\eta)$, is a continuous function of $\eta$ which is zero unless
$\eta(x)=1-\eta(y)=1$.  

Let
 \begin{equation}
  \label{5}
  \X_n = \{\eta\in \{0,1\}^Z: \sum_{x\le n} \eta(x) 
      = \sum_{x>n} (1 - \eta(x))<\infty\}
 \end{equation}
and 
 \begin{equation}
  \label{5a}
  \X= \cup_n \X_n. 
 \end{equation}
 The set $\X$ is countable; we will call elements of $\X$ {\it blocking
configurations}, and call probability measures supported on $\X$ {\it
blocking measures}.  Our interest is to find sufficient conditions on the
rates $p$ for the existence of blocking measures which are invariant for the
process $\eta_t$. 

We will construct the process on blocking configurations directly.  For the
construction we will use two conditions on the rates, which we assume
throughout the paper; these could be somewhat weakened at the price of
increasing the complexity of the exposition.  (Liggett \cite{[Ligg]} gives
conditions on the rates which assure the existence of the process started
from an arbitrary initial condition.) First, we take the rates to be
uniformly bounded; we can set the upper bound equal to one by a time-scale
change, and thus assume that
 \begin{equation}
\label{68}
0\le p(x,y;\eta)\le1.
 \end{equation}
 Second, we assume that the total rate for exiting any configuration of $\X$
is finite: 
 \begin{equation}
 \label{666}
  \hbox{For any }\zeta\in\X,\;\;\;\;
      \sum_x \sum_y \bar p(x,y;\zeta)\,<\,\infty. 
 \end{equation}
 This condition follows from \reff{68} if there is an upper bound on the
range of jumps. 

Various special cases are of interest.  The rates are {\it simple} when they
are independent of the configuration except for the exclusion condition, so
that
 \begin{equation}
 \label{S}
   p(x,y;\eta)=c(x,y)\eta(x)(1-\eta(y)),
 \end{equation}
  and are {\it translation invariant} when
 \begin{equation}\label{TI}
   p(x,y;\eta)=p(0,y-x;\tau_{-x}\eta),
 \end{equation} 
 where $\tau_z$ is the operator of translation by $z$.  When both of these
conditions are satisfied the rates can be written in the form 
 \begin{equation}\label{STI}
   p(x,y;\eta)=a(y-x)\eta(x)(1-\eta(y)).
 \end{equation} 
  Liggett \cite{[Ligg]} exhibits invariant blocking measures 
in the case of simple,
translation invariant rates with jumps restricted to length 1: $a(z)=0$
for $|z|>1$ and $a(1)>a(-1)$.  A trivial extension of his result is the
following: if for some $\alpha<1$ the rates have the form \reff{S} with
 \begin{equation}
  \label{3}
  c(x,y) = \alpha^{x-y}\, c(y,x) \qquad\hbox{for all $x<y$},
 \end{equation}
then the product measure $\mu$ with marginals
 \begin{equation}
  \label{4}
  \mu(\eta(x) = 1)= {1 \over 1+\alpha^x}
 \end{equation}
 is reversible for the process $\eta_t$.  This is a special case of a more
general construction which we describe in the appendix. 
 
For a more general set of rates, one might expect that blocking
measures exist when the process has a
sufficiently strong positive drift, for example in the translation
invariant simple case (that is, for rates satisfying \reff{STI}) when
 \begin{equation}
  \label{2}
  \sum_{\{y\}} \,y\,a(y)\, >0\,,
 \end{equation}
 (positive mean drift for the underlying random walk).  Proving that
(\ref{2}) or a similar condition implies the existence of blocking measures
seems quite difficult; this is one of the open problems of \cite{[Ligg]}.
When the rates $p(x,y;\eta)$ depend on the configuration $\eta$ at sites
other than $x$ and $y$, it is not even clear what necessary and/or sufficient
condition to conjecture.  We do not deal directly with conditions like
\reff{2}, but give a different sort of sufficient condition, showing that
when the rates of two processes are appropriately related, existence of a
blocking measure for one implies existence for the second. 

Note that if $\mu$ is any invariant blocking measure for $\eta_t$ then
$\mu(\X_n)\ne0$ for some $n$; since each $\X_n$ is a closed set for the
process, the conditional measure $\mu_n=\mu(\cdot|\X_n)$ is then also an
invariant blocking measure.  Thus, if we permit ourselves a translation of
the entire system, there is no loss of generality in treating existence of a
blocking measure on $\X$ as equivalent to the existence of a blocking measure
on $\X_0$.  We remark that if the rates are simple and translation invariant
(see \reff{STI}) then $\X_n$ is irreducible whenever there is a positive rate
for some forward and some backward jump, and the greatest common divisor of
$\{\,x\ne0\,:\,a(x)>0\,\}$ is 1, so that under these condition each $\mu_n$
is unique and extremal in the class of invariant blocking measures. 

We now compare the process $\eta_t$ with a second process $\bar\eta_t$ for
which the generator $\bar L$ is constructed as in (\ref 1) but with rates
$\bar p(x,y;\eta)$.  Our main result, presented in
Section~\ref{Result}, gives conditions on the rates $p$ and $\bar p$
under which the existence of a blocking invariant measure for the process
$\bar\eta_t$ implies the existence of such a measure for $\eta_t$.  In the
case in which the rates are simple and translation invariant, it takes the
following form:

 \begin{thm}\label{6} Suppose that
$p(x,y;\eta)=a(y-x)\eta(x)(1-\eta(y))$ and that
 \begin{eqnarray}
  \label{7}
 a(x) \ge \bar a(y)\,, &\ \ \ \hbox{for}&  0<x\le y\,, \\
 a(y) \le \bar a(x)\,, &\ \ \ \hbox{for}&  y\le x<0\,. \label{7a}
 \end{eqnarray}
 Then if $\bar\eta_t$ has a blocking invariant measure, so does $\eta_t$. 
 \end{thm}

For example, we may take the weights $\bar p$ to have the form \reff{3},
with $c(x,y)=a(y-x)$ for $x<y$, as in \reff{STI}, so that the requisite
blocking measure is given by \reff{4}.  

We remark that establishing the existence of invariant blocking measures is a
special case, and perhaps a first step toward the general case, of the
problem of establishing the existence of {\it invariant shock measures}:
measures on $\{0,1\}^\Z$ which have distinct asymptotic limits to the right
and left of the origin and which are time invariant in some appropriate
sense, usually for the process as seen from a suitable random viewpoint. 
Such measures are related to the shock solutions of the Burgers equation,
which describes the process in the hydrodynamical limit.  The left and right
asymptotic measures will be time invariant for the process in the usual
sense, so that invariant shock measures appear in systems that have more than
one translation invariant state.  Given two such asymptotic measures, the
shock measure describes one ultimate fate of the system when it starts with
one of these on each side of the origin (another is the so called {\it
rarefaction fan}).  The blocking measures are the simplest shock measures:
conceptually, because they are invariant when seen from a fixed viewpoint,
and technically, because they have support on a countable state space. 

In the case of simple exclusion the extremal time and translation invariant
measures are the one parameter family of homogeneous product measures indexed
by density.  In nearest neighbor asymmetric simple exclusion, existence of
invariant shock measures has been established for the process as seen from a
``second class particle'', (\cite{[FKS],[F],[DJLS],[DLS]}).  The approach of
\cite{[FKS]} and \cite{[F]} was closely based on the known blocking measures
for this process, the product measures \reff{4}.  In \cite{[DGLS]} other
approaches for the problem of describing shock measures are proposed. 

The paper is organized as follows.  In Section~\ref{Construction} we
construct $\eta_t$ on ${\cal X}$ using Poisson processes (the Harris
graphical construction); the construction is made in
such a way as to facilitate an appropriate coupling of two such
process.  We describe in Section~\ref{Order} the key idea for the
proof of our results: the introduction of a certain partial
order~$\prec$ on the space $\X_0$ of blocking configurations with the
property that, under the coupling, the conditions of Theorem~\ref{6}
(or the more general conditions to be given later) imply that if the
initial configurations $\eta_0$ and $\bar \eta_0$ satisfy
$\eta_0\prec\bar\eta_0$, then this ordering is preserved by the
dynamics: $\eta_t\prec\bar\eta_t$ for all $t\ge0$.  In
Section~\ref{Result} we state and prove our general result, of which
Theorem~\ref{6} is an immediate corollary.  In
Section~\ref{Applications} we give some applications, and in the
appendix discuss the construction of a class of possible comparison
processes $\bar\eta_t$.

\section{Construction of the process}
\label{Construction}

We exhibit now a special construction of the process in ${\cal X}$.
The construction requires that the rates $p(x,y;\zeta)$ satisfy conditions
\reff{68} and \reff{666} of the introduction.

For a configuration $\eta\in \X$ we define ordered positions of the
particles and empty sites by
 \begin{eqnarray}
  \label{8}
 && x_0(\eta) = \min\{x:\eta(x) =1\}\,,\\
&& x_k(\eta) = \min\{x>x_{k-1}(\eta):\eta(x) =1\}\,,\\
 && y_0(\eta) = \max\{x:\eta(x) =0\}\,,\\
&& y_k(\eta) = \max\{x<y_{k-1}(\eta):\eta(x) =0\}\,.\label{8bis}
 \end{eqnarray}
 For each pair $(i,j)$ with $i,j\geq 0$ let
 \begin{equation}
  \label{20}
 \Theta^{i,j}\,=:\,
 \{((T_n^{i,j},U^{i,j}_n), (R_m^{i,j},V^{i,j}_m)):n,m\ge 1\}
 \end{equation}
 be a process with the following properties:
\begin{itemize}
\item Both $(T_n^{i,j}-T_{n-1}^{i,j})_{n\ge1}$ and
$(R_m^{i,j}-R_{m-1}^{i,j})_{m\ge1}$, where by convention
$T_0^{i,j}=R_0^{i,j}=0$, are families of independent exponentially
distributed random variables of mean one. In other words, $(T_n^{i,j})$ and
$(R_m^{i,j})$ are Poisson processes of rate $1$ for all $i,j$.
\item Both $(U_n^{i,j})_{n\ge1}$ and
$(V^{i,j}_m)_{m\ge1}$ 
are families of independent random variables, uniformly distributed in
$[0,1]$.
\item All four of these families of variables are mutually
independent. 
\end{itemize}
We also assume that  $\{\Theta^{i,j}:i,j\ge 0\}$ is a family of mutually
independent processes.
The times $T^{i,j}_n$ and $R^{i,j}_m$ will be called Poisson events and the
associated random variables $U^{i,j}_n$ and $V^{i,j}_m$ will be called
{\sl marks}. 

We now construct the process $\eta_t$ as a function of
the marked Poisson processes and the initial configuration $\eta_0\in\X$.  
Set $\tau_0=0$ and suppose inductively that we have
defined times $\tau_0,\ldots,\tau_{n-1}$ and configurations
$\eta_{\tau_0},\ldots,\eta_{\tau_{n-1}}$.
Define 
 \begin{equation}
\label{taun}
  \tau_n= \min\Bigl\{\inf_{i,j,k} \{T_k^{i,j}>\tau_{n-1}: U_k^{i,j}<
     A_+(\eta_{\tau_{n-1}},i,j)\},\,
   \inf_{i,j,k} \{R_k^{i,j}>\tau_{n-1}: V_k^{i,j}<
     A_-(\eta_{\tau_{n-1}},i,j)\}\Bigr\}\,, \label{17}
 \end{equation}
where for $i,j\ge0$,
 \begin{eqnarray}
  A_+(\eta,i,j) &=& p(x_i(\eta),y_j(\eta);\eta)\,
       \one\{y_j(\eta)>x_i(\eta)\},\label{16}\\
  A_-(\eta,i,j) &=& p(x_i(\eta),y_j(\eta);\eta)\,
       \one\{y_j(\eta)<x_i(\eta)\}.\label{16bis}
  \end{eqnarray}
Here $\one S$ denotes the characteristic function of the set $S$. 
If $(I_n,J_n)$ is the pair $(i,j)$ such that $T_k^{i,j}$ or
$R_k^{i,j}$ realizes the infimum $\tau_n$ for some $k$, set
\begin{eqnarray}
  \label{eq:xn}
  X_n&=&x_{I_n}(\eta_{\tau_{n-1}}),\\
 Y_n&=&y_{J_n}(\eta_{\tau_{n-1}}),
\end{eqnarray}
and define
 \begin{eqnarray}
  \eta_{\tau_n}
  &=& (\eta_{\tau_{n-1}})^{X_n,Y_n}.
 \end{eqnarray}
 This completes the induction step.  To finish the construction after all 
$\tau_n$ and $\eta_{\tau_n}$ are defined, set
 \begin{equation}
  \label{11}
  \eta_t =  \sum_{n\ge 0}\eta_{\tau_n}\one\{\tau_n\le t<\tau_{n+1}\}
  \qquad\hbox{for all $t\ge0$}.
 \end{equation}
 It is important to notice that after each jump the particles and holes are
effectively relabeled according to \reff{8}--\reff{8bis}, so that for all
times $t$,
 \begin{equation}
  \label{eq:400}
  x_i(\eta_t) \le x_{i+1}(\eta_t) \quad\hbox{ and }\quad
    y_j(\eta_t) \ge y_{j+1}(\eta_t),\qquad i,j\ge0.
 \end{equation}

The construction may be described in words as follows.  We use independent
times ($T^{i,j}_n$ and $R^{i,j}_m$, respectively) for jumps to the right and
jumps to the left; this is not necessary for the construction here but
ensures that the coupling we define later preserves a certain partial order
on configurations.  The instant $\tau_n$ is the first time after $\tau_{n-1}$
at which a jump is performed, and is the minimum of the first scheduled jump
times to the right and to the left.  The first scheduled jump time to the
right is the first $T^{i,j}_k$ for which the corresponding uniform random
variable $U^{i,j}_k$ is smaller than the threshold $A_+$, defined by
\reff{16} to ensure that the jump is indeed to the right and occurs at the
correct rate (here we use the condition \reff{68} that $p(x,y;\eta)\le1$). 
Similarly, the first scheduled jump time to the left is the first $R^{i,j}_k$
for which the corresponding uniform random variable $V^{i,j}_k$ is smaller
than the threshold $A_-$ defined by \reff{16bis}.  The configuration at time
$\tau_n$ is then the one obtained by interchanging the hole and the particle
whose indexes $i,j$ correspond to the $R^{i,j}_k$ or $T^{i,j}_k$ that
realizes the time $\tau_n$. 

To see that the above is well defined for initial configurations in $\X$ it
suffices to see that, for any initial $\eta_0\in \X$, $\tau_n$ is with
probability one a strictly increasing sequence of (finite) times.  The
conditional distribution of $\tau_n-\tau_{n-1}$ given the past up to
$\tau_{n-1}$ is
 \begin{eqnarray}
\label{300}
\P(\tau_n - \tau_{n-1} > s\,\vert\, \eta_{\tau_{n-1}})
 &=&
\exp\Bigl\{-\sum_{x,y} p(x,y;\eta_{\tau_{n-1}})\Bigr\},
 \end{eqnarray}
 by \reff{taun} (it is the minimum of independent random variables with
 exponential distribution and inverse-mean $p(x,y;\eta_{\tau_{n-1}})$).
 Since $\eta_{\tau_{n-1}}$ is obtained by doing at most $n-1$
 modifications to the initial configuration $\eta_0$, it belongs to
 $\X$.  By condition \reff{666} of the introduction, the conditional law
\reff{300}
 is that of a non-degenerate exponential random variable.
 
 It is tedious but easy to show that the process $\eta_t$ so
 constructed in $\X$ has generator $L$ (restricted to $\X$). We
 remark that the above construction works also if the process
 restricted to $\X$ has explosions, that is, if $\lim_{n\to\infty}
 \tau_n<\infty$.

We give now a graphical interpretation of this construction, and of
the coupling of the processes to be introduced later. For simplicity
assume $\X=\X_0$. To each configuration $\eta\in\X_0$ associate an
interface $\Phi\eta$ corresponding to the integrated profile of
$\eta$.  Here $\Phi:\X_0\to\Z_+^\Z$ is defined by either of two
equivalent expressions:
 \begin{eqnarray}
  \label{18}
  (\Phi\eta)(x) &=& -x + 2 \sum_{y\le x}\eta(y)\\
&=& x + 2 \sum_{y>x}(1-\eta(y))\label{18a}
 \end{eqnarray}
 Note that $\Phi\eta$ increases by one when a particle is present at
 $x$ or decreases by one when no particle is present at $x$, so that
 in particular, $\vert\Phi\eta(x)-\Phi\eta(x+1)\vert = 1$.  The graph
 $\{\,(x,(\Phi\eta)(x))\mid x\in\Z\,\}$ is a subset of the lattice
 $\Z^2_{\rm even}=\{\,(x,y)\in\Z^2\mid x+y\hbox{ is even}\,\}$.  The
 Heaviside configuration $\eta^H$, given by $\eta^H(x) = \one\{x\ge
 1\}$, gives rise to the interface $\Phi\eta^H(x) = |x|$.

The interface picture yields a geometric interpretation of the construction
of the process $\eta_t$.  Index the squares (plaquettes) of the lattice
$\Z^2_{\rm even}$ as $\{\,S_{i,j}\mid i,j\in\Z\,\}$ as shown in Figure~1
($S_{i,j}=\{\,(x,y)\mid2i<x+y<2i+2, 2j<y-x<2j+2\,\}$); with this convention,
the interface $\Phi\eta$ lies above $S_{i,j}$ ($i,j\ge0$) if and only if
$x_i(\eta)<y_j(\eta)$.  Now think of the marked processes
$(T_n^{i,j},U^{i,j}_n)$ and $(R_m^{i,j},V^{i,j}_m)$ as associated with
$S_{i,j}$.  When at the Poisson event $T^{i,j}_n$ the
corresponding uniform variable $U^{i,j}_n$ is less than
$p(x_i(\eta),y_j(\eta);\eta)$, then, if the interface $\Phi\eta$ lies above
$S_{i,j}$, we update the interface by decreasing its height by two units in the
interval $(x_i,y_j]$.  Similarly, when at time $R^{i,j}_m$ the corresponding
mark satisfies 
$V^{i,j}_m<p(x_i(\eta),y_j(\eta);\eta)$ and the interface lies below
$S_{i,j}$, we increase by two units the height of the interface in
the interval $(x_i,y_j]$.  All of this is shown in Figure~1.

\section{An order relation on configurations}
\label{Order}

 For configurations $\eta$ and $\bar\eta\in\X_0$, we say that
 \begin{equation}
  \label{eq:30}
 \eta\prec\bar\eta\ \  \hbox{ if and only if \ \ for all $i,j\ge0$, }
x_i(\eta)\ge x_i(\bar\eta)\hbox{ and }y_j(\eta)\le y_j(\bar\eta).
 \end{equation}
 It is easy to see that this is a partial order which corresponds to the
natural order on interfaces:
 \begin{equation}
  \label{19}
  \eta\prec\bar\eta \ \ \ \hbox{ if and only if }\ \ \  (\Phi\eta)(x) \le
  (\Phi\bar\eta)(x) \hbox{ for all }x\in\Z.
 \end{equation}
 Under this ordering, the Heaviside configuration $\eta^H$ precedes every
other configuration: $\eta^H\prec\eta$ for any $\eta\in\X_0$.  From \reff{18}
and \reff{18a} it follows that if $\eta\prec\bar\eta$ then for all $z\in\Z$,
 \begin{eqnarray}
  \label{19f} (\Phi\bar\eta)(z) - (\Phi\eta)(z)&=& 2\,\sum_i \one\{
  x_i(\bar\eta) \le z < x_i(\eta)\}\\ &=& 2\,\sum_j \one\{
  y_j(\eta) \le z < y_j(\bar\eta)\}\,,\label{19g}
 \end{eqnarray}
 and for all $x,y$ such that $\eta(x)=1$ and $\eta(y)=0$ and all $z\in\Z$,
 \begin{equation}
  \label{21}
  (\Phi\eta^{x,y})(z) = (\Phi\eta)(z) - 2\,\one\{x \le z < y \} +
  2\,\one\{y \le z < x \}.
 \end{equation}

 The following lemma says essentially that if we have two configurations which
are ordered by $\prec$ then they will remain ordered after either (i)~a jump
in both configurations, in the same direction, of the $i^{\rm th}$ particle
to the $j^{\rm th}$ hole, or (ii)~certain jumps in only one of the
configurations. 

 \begin{lemma}
  \label{32} Assume $\eta\prec\bar\eta$, fix $i$ and $j$, and let
$x=x_i(\eta)$, $y=y_j(\eta)$, $\bar x= x_i(\bar\eta)$, and $\bar
y=y_j(\bar\eta)$.  Then jumps preserve ordering in the following cases:
 \begin{eqnarray}
\label{340} 
&&\hbox{If $\bar x \le x < y \le \bar y$, then
$\eta^{x,y}\prec\bar\eta$.} \\
\label{341} 
&&\hbox{If $y\le \bar y< \bar x\le x $, then $ \eta\prec
\bar\eta^{\bar x,\bar y}$}. \\
\label{240} 
&&\hbox{If $\bar x \le x < y \le \bar y$, then
$\eta^{x,y}\prec\bar\eta^{\bar x,\bar y}$ } \\
\label{241} 
&&\hbox{If $y\le \bar y< \bar x\le x $, then $\eta^{x,y}\prec
\bar\eta^{\bar x,\bar y}$.} \\
\label{242}
&&\hbox{If $x>y$ and $\bar x<\bar y$, then $\eta\prec \bar\eta^{\bar
x,\bar y}$ and $ \eta^{x,y}\prec \bar\eta $}.
 \end{eqnarray}
 \end{lemma}

Before giving a formal proof of this lemma, we describe its graphical
interpretation.  The interface $\Phi\eta$ lies below $\Phi\bar\eta$.  In
cases \reff{340} and \reff{240} the square $S_{i,j}$ lies below both
interfaces, so that for either interface a jump of the $i^{\rm th}$
particle to the $j^{\rm th}$ hole---briefly, an $(i,j)$ jump---lowers the
interface; \reff{340} and \reff{240} assert respectively that the order is
preserved by either a jump in the lower interface only, or a jump for both
interfaces.  Similarly, in cases \reff{341} and \reff{241} $S_{i,j}$ lies
above both interfaces, an $(i,j)$ jump raises either interface, and the
order is preserved by such a jump in either the upper interface alone or in
both.  Finally, in case \reff{242} $S_{i,j}$ lies between the two
interfaces,  an $(i,j)$ jump for the lower interface raises it and for the
upper interface lowers it, and \reff{242} asserts that such a jump for
either interface alone preserves the order.  These properties are easy to
check in the graphical representation.

 \vskip10pt
 \noindent
 {\bf Proof of Lemma~\ref{32}:} Statements \reff{340} and \reff{341}
follow immediately from \reff{21}.  Under the hypothesis of \reff{240} $x<y$
and $\bar x <\bar y$.  Hence, by \reff{21},
 \begin{eqnarray}
  \label{243} (\Phi\eta^{x,y})(z) 
    &=& (\Phi\eta)(z) - 2\,\one\{x \le z < y \};
 \end{eqnarray}
an analogous identity holds for $\bar\eta$.  Since
$\eta\prec\bar\eta$ and $\bar x\le x<y\le \bar y$, by \reff{19f} and \reff{19g},
 \begin{eqnarray}
&&(\Phi\eta)(z) \le (\Phi\bar\eta)(z)
    - 2\,\one\{\bar x \le z < x\}-2\,\one\{y \le z < \bar y\}.
 \end{eqnarray}
Subtracting $2\,\one\{x<z\le y\}$ in both members of the above
inequality we get
 \begin{eqnarray}
  \label{244} (\Phi\eta)(z) - 2\,\one\{x \le z < y \}&\le&
(\Phi\bar\eta)(z) - 2\,\one\{\bar x \le z < \bar y \},
 \end{eqnarray}
 which by \reff{21} is the same as 
$(\Phi\eta)(z)<(\Phi\bar\eta^{\bar x,\bar y})(z)$.
In this way we get
$\eta^{x,y}\prec\bar\eta^{\bar x,\bar y}$ and \reff{240} is proven.
Display \reff{241} is verified analogously.

By \reff{19f}, 
 \begin{eqnarray}
\label{350}
 (\Phi\bar\eta)(z) - (\Phi\eta)(z) &\ge& 2\,\one\{y \le z < \bar y\}\,,\\
(\Phi\bar\eta)(z) - (\Phi\eta)(z) &\ge& 2\,\one\{\bar x \le z < x\}\,.
 \end{eqnarray}
Under the hypothesis of \reff{242}, this implies that
 \begin{eqnarray}
\label{351}
 (\Phi\bar\eta)(z) - (\Phi\eta)(z)
   &\ge& 2\,\one\{\min\{y,\bar x\} \le z < \min\{\bar y, x\}\}\,.
 \end{eqnarray}
 Applying \reff{21}, we get \reff{242}. \square

\section{Statement and proof of main result}
\label{Result}

 Now we consider two processes $\eta_t$ and $\bar\eta_t$ with rates $p$ and
$\bar p$, respectively, as discussed in the introduction.  Our main result
is:

 \begin{thm}\label{66} Suppose that whenever $\eta\prec\bar\eta$ and
$\eta(x)=\bar\eta(\bar x)=1$, $\eta(y)=\bar\eta(\bar y)=0$, 
 \begin{eqnarray}
  \label{70}
 p(x,y;\eta) \ge \bar p(\bar x,\bar y;\bar\eta)\,,&
         \ \ \ \hbox{if}&  \bar x\le x<y\le\bar y\,, \\
 p(x,y;\eta) \le \bar p(\bar x,\bar y;\bar\eta)\,,&
         \ \ \ \hbox{if}&  y\le \bar y < \bar x\le x\,.\label{70a}
 \end{eqnarray}
 Then if $\bar\eta_t$ restricted to $\X$ has a blocking invariant
 measure, so does $\eta_t$.
 \end{thm}

 Theorem~\ref{6} is an immediate corollary of Theorem~\ref{66}.
 
 We construct simultaneously the two processes $\eta_t$ and
 $\bar\eta_t$ using the {\it same} marked Poisson processes
 $((T^{i,j}_n,U^{i,j}_n),(R^{i,j}_m,V^{i,j}_m))$.  This joint
 construction is called {\em coupling} and is the key to the proof.

 \begin{lemma}
\label{12} 
Assume that $\eta_t$ and $\bar\eta_t$ are processes with rates $p$ and
$\bar p$ satisfying \reff{70}---\reff{70a}. Under the coupling, if
$\eta_0\prec \bar\eta_0$ are both configurations of $\X$, then for
all $t\ge 0$, $\eta_t\prec \bar\eta_t$.
 \end{lemma}
 
 \proof This is a mark-by-mark proof.  Set $\theta_0=0$ and let
 $\theta_1<\theta_2<\cdots$ be the instants at which there is a jump
 for at least one of the processes $\eta_t,\bar\eta_t$.  Assume
 inductively that $\eta_{\theta_{n-1}}\prec\bar\eta_{\theta_{n-1}}$,
 so that if $(x_i,y_j)$ and $(\bar x_i,\bar y_j)$ are the sites and
 holes of $\eta_{\theta_{n-1}}$ and $\bar\eta_{\theta_{n-1}}$,
 respectively, at time $\theta_{n-1}$, then
 \begin{equation}
  \label{13}
  x_i\ge \bar x_i,\quad\hbox{and}\quad  y_j\le \bar y_j,\qquad i,j\ge0.
 \end{equation}
Let $\tau_n$ and $\bar\tau_n$ be the times defined as in \reff{17} for the
processes $\eta_t$ and $\bar\eta_t$, so that
 \begin{equation}
  \label{14}
  \theta_n = \min \Bigl\{ \min \{\tau_k>\theta_{n-1}\}, \min
  \{\bar\tau_k>\theta_{n-1}\}\Bigr\}.
 \end{equation}
 Let $(I,J,K)$ be the indices which realize the infimum \reff{17} defining
the time $\theta_n$, so that $\theta_n \in\{T_K^{I,J},R_K^{I,J}\}$.  Let
$U\in\{ U_K^{I,J},V_K^{I,J}\}$ be the uniform random variable related with
the indexes realizing the infimum, and let $\sigma=\pm$ indicate the
direction of the jump at $\theta_n$: $\sigma=+$ if $\theta_n=T_K^{I,J}$ and
$U=U_K^{I,J}$, $\sigma=-$ if $\theta_n=R_K^{I,J}$ and $U=V_K^{I,J}$.  Let
 \begin{eqnarray}
  &&X = x_I,\,\,\,\, \bar X = \bar x_I,\,\,\,\, Y= y_J ,\,\,\,\,
    \bar Y=\bar y_J;\\ 
  &&\xi =\eta_{\theta_{n-1}} \;\; \; \; \; \bar\xi=\bar\eta_{\theta_{n-1}};
   \label{15}\\
  && B= A_\sigma(\xi,I,J),\;\; 
    \bar B = A_\sigma(\bar\xi,I,J).
 \end{eqnarray}

Since \reff{13} implies that 
$\bar X \le X$ and $Y\le \bar Y$, there are three possibilities:

 \begin{enumerate}
\item $\bar X \le X < Y \le \bar Y$. By hypothesis \reff{70}, $\bar B\le B$.
  Hence there are two possibilities:

\subitem{(a)} $U < \bar B \le B$. In this case $\eta_{\theta_n} =
  \xi^{X,Y}$ and $\bar\eta_{\theta_n} = \bar\xi^{\bar X,\bar Y}$. By
  \reff{240}, $ \eta_{\theta_n}\prec \bar\eta_{\theta_n}$.

\subitem{(b)} $\bar B \le U < B$. In this case $\eta_{\theta_n} =
  \xi^{X,Y}$ and $\bar\eta_{\theta_n} = \bar\xi$. By
  \reff{340}, $ \eta_{\theta_n}\prec \bar\eta_{\theta_n}$.

\item $Y\le \bar Y< \bar X\le X $. By hypothesis \reff{70a}, $B\le \bar B$.
  Hence there are two possibilities:

\subitem{(a)} $U < B \le \bar B$. In this case $\eta_{\theta_n} =
  \xi^{X,Y}$ and $\bar\eta_{\theta_n} = \bar\xi^{\bar X,\bar Y}$. By
  \reff{241}, $ \eta_{\theta_n}\prec \bar\eta_{\theta_n}$.

\subitem{(b)} $B \le U < \bar B$. In this case $\eta_{\theta_n} =
  \xi^{X,Y}$ and $\bar\eta_{\theta_n} = \bar\xi$. By
  \reff{341}, $ \eta_{\theta_n}\prec \bar\eta_{\theta_n}$.

\item $X>Y$ and $\bar X<\bar Y$. There are two possibilities: 

\subitem{(a)} $\sigma=+$ and $0=B\le U=U^{I,J}_K < \bar B$.
  In this case $\eta_{\theta_n} = \xi$ and 
  $\bar\eta_{\theta_n} = \bar\xi^{\bar X,\bar Y}$. By
  \reff{242}, $ \eta_{\theta_n}\prec \bar\eta_{\theta_n}$.

\subitem{(b)} $\sigma=-$ and $0=\bar B\le U=V^{I,J}_K < B$.
  In this case $\eta_{\theta_n} = \xi^{X,Y}$ and 
  $\bar\eta_{\theta_n} = \bar\xi$. Again by
  \reff{242}, $ \eta_{\theta_n}\prec \bar\eta_{\theta_n}$. \square
 \end{enumerate}

Notice that if we had used the same Poisson process for both forward
and backward jumps then in the situation of case 3 above jumps could
have occurred simultaneously in $\eta$ and $\bar\eta$, in opposite
directions, which could destroy the ordering.

We remark that explosions are not excluded in Lemma \ref{12}.

\vskip 3mm

 \noindent
 {\bf Proof of Theorem~\ref{66}:} As remarked in the introduction, it
suffices to show that if $\bar\eta_t$ has an invariant measure in $\X_0$, so
does $\eta_t$.  By restricting to a subset of $\X'\subset \X_0$ (if
necessary) we may assume that $\bar\eta_t$ is ergodic with invariant measure
$\bar\mu$ having support $\X'$.  This excludes explosions for the process
$\bar\eta_t$ starting with configurations in $\X'$. 

Start the coupled process with any two configurations $\zeta \prec \bar \zeta$,
with $\bar\zeta\in\X'$ and $\zeta\in\X$. We know that:
\begin{enumerate}
\item $\eta_t \prec \bar \eta_t$, by Lemma~\ref{12};
\item No explosions occur for $\eta_t$ (by an argument similar to the
one in the proof of Lemma \ref{12});
\item Since $\bar\eta_t$ is a continuous time ergodic Markov process in a
  countable state space, it converges in distribution to its unique invariant
  measure $\bar \mu$.
\end{enumerate}
Hence any weak Cesaro-limit $\mu$ of the distribution of $\eta_t$ is coupled
with $\bar \mu$ in such a way that, calling $\nu$ the coupled measure with
marginals $\mu$ and $\bar\mu$, $\nu$ satisfies
 \begin{equation}
   \label{eq:402}
   \nu((\eta,\bar\eta): \eta\prec\bar\eta) =1\,.
 \end{equation}
 This in particular implies that $\mu(\X)=1$.  Since $\mu$ is a
Cesaro-limit, $\mu$ is invariant for $\eta_t$.  This implies the theorem. 
\square

\section{Applications}
\label{Applications}

To apply Theorem~\ref{66} one needs a suitable comparison process $\bar\eta$
which is known to have an invariant blocking measure.  Obvious candidates are
processes satisfying \reff{3}, for which the product measures \reff{4} are
invariant; in this section we draw some simple conclusions from this
comparison.  In the appendix we  discuss briefly the existence of
other possible comparison processes: those which satisfy detailed balance
with respect to a Gibbs measure obtained from a suitable potential
(Hamiltonian).

 \begin{thm} 
 \label{app1} 
 Suppose that the exclusion process $\eta_t$ has simple, translation
invariant rates $p(x,y;\eta)=a(y-x)\eta(x)(1-\eta(y))$ which for some
$\alpha$ with $0\le\alpha<1$ satisfy
 \begin{equation}
  a(-x)\le\alpha^x\displaystyle\inf_{0<y\le x}a(y)
 \end{equation}
  for all $x>0$.  Then
$\eta_t$ has an invariant blocking measure. 
 \end{thm}

\proof The process with rates
 $\bar p(x,y;\eta)=\bar a(y-x)\eta(x)(1-\eta(y))$, where for $x>0$,
 \begin{equation}
\bar a(x)=\inf_{0<y\le x}a(y)\qquad\hbox{and}\qquad
   \bar a(-x)= \alpha^x\bar a(x)\,,
 \end{equation}
 has an invariant measure of the form \reff{4}.  Thus the process $\eta_t$
has an invariant blocking measure by Theorem~\ref{6}. \square
 
As a second example, consider a process with symmetric ``disorder,'' in which
translation invariant, asymmetric, nearest neighbor rates are perturbed by
arbitrary, bounded, symmetric nearest neighbor rates. Specifically, take 
$p(x,y;\eta)=(c_0(x,y)+c_1(x,y))\eta(x)(1-\eta(y))$, where
$c_1(x,y)=c_2(x,y)=0$ if $|x-y|>1$ and
 \begin{equation}
   c_0(x,x+1) = K, \qquad c_0(x+1,x)=0,\qquad
  c_1(x,x+1) = c_1(x+1,x)=h(x), 
 \end{equation}
 with $K>0$ and $h:\Z\to\R_+$ an arbitrary bounded function.  It follows
from Theorem~\ref{66} that this process has a blocking measure.  A suitable
comparison process has rates $\bar p(x,y;\eta)=\bar c(x,y)\eta(x)(1-\eta(y))$
with $\bar c(x,x+1)=c(x,x+1)$, $\bar c(x+1,x)=\alpha c(x,x+1)$, and $\bar
c(x,y)=0$ if $|x-y|>1$, where $\alpha=M/(M+K)$ with $M$ an upper bound
on $h(x)$; these rates satisfy \reff{3} and hence have a blocking measure
as given in \reff{4}.  We single out this rather trivial example because in
this case it is easy to see that the product measures with constant density
are invariant measures, since if $\mu$ is such a measure then
$L_1^*\mu=L_2^*\mu=0$ and hence $L^*\mu=0$, where $L_i^*$ is the adjoint of
the generator for the process with rates $c_i$.

\section*{Acknowledgments} JLL was supported in part by NSF Grant
DMR-9813268.  PAF and JLL thank DIMACS and its supporting agencies,
the NSF under contract STC-91-19999 and the N. J. Commission on Science and
Technology. This work started while PAF. was visiting Rutgers University
with support of DIMACS. PAF was supported in part by FAPESP.

\section*{Appendix}

\setcounter{equation}{0}
\setcounter{thm}{0}
\def\theequation{A.\arabic{equation}}
\def\thethm{A.\arabic{thm}}

The remark that processes satisfying \reff{3} have invariant product blocking
measures of the form \reff{4} can be generalized to processes which satisfy
detailed balance with respect to a Gibbs measure obtained from a suitable
potential (Hamiltonian).  The latter is specified \cite{[Ligg]} by a
collection of real numbers $\{J_R\}$ indexed by finite subsets $R$ of $\Z$
and satisfying $\sum_{R\ni x}|J_R|<\infty$ for each $x\in\Z$.  We show that
if these coupling constants are chosen appropriately, then blocking Gibbs
measures for this potential arise as the limit of finite volume measures. 

Let $T_N=[-N+1,N]\cap\Z$ and $Y_N=\{0,1\}^{T_N}$.  For
$\eta\in Y_N$ let $\eta^*\in\X$ be the configuration which agrees
with $\eta$ in $T_N$ and with $\eta^H$ outside $T_N$.  The {\it energy} of
the configuration $\eta$ is
 \begin{equation}
H_N(\eta)=\sum_{\{R\mid R\cap T_N\ne\emptyset\}}J_R\;\chi_R(\eta^*),
 \end{equation}
 where $\chi_R(\zeta)=\prod_{x\in R}(2\zeta(x)-1)$; the variables
$2\zeta(x)-1$ are spins which take values $\pm1$.  The corresponding
finite-volume {\it Gibbs measure} $\nu_N$ on $Y_N$ is defined
by
 \begin{equation}
\nu_N(\{\eta\}) =Z_N^{-1}\exp\bigl(-H_N(\eta)\bigr)
 \end{equation}
 for $\eta\in Y_N$, with $Z_N=\sum_{\zeta\in Y_N}\exp(-H(\zeta))$ a
normalization constant; $\nu_N$ defines a measure on $\{0,1\}^\Z$ by
setting $\nu_N(A)=\nu_N(\{\eta\in\{0,1\}^\Z\mid \eta^*\in A\}$.  Now
let us assume for simplicity that all $n$-body terms in the potential, for
$n\ge2$, are translation invariant, i.e., that $J_{R+k}=J_R$ for $k\in\Z$ and
$|R|\ge2$ (this assumption could easily be relaxed), and let $K=\sum_{R\ni
x,\,|R|\ge2}|J_R|$.  

 \begin{thm}\label{tdlim} Suppose that the one particle potential $J_{\{x\}}$
approaches $\mp\infty$ as $x$ approaches $\pm\infty$, respectively,
sufficiently fast that
 \begin{equation}
  \label{conv}
   \sum_{x\ge1} \exp\bigl(2J_{\{x\}}\bigr)<\infty \qquad\hbox{and}\qquad
   \sum_{x\le0} \exp\bigl(-2J_{\{x\}}\bigr)<\infty.
 \end{equation}
 Then $\nu=\lim_{N\to\infty}\nu_N$ exists and is a blocking measure. 
Moreover, if the rates $p(x,y;\eta)$ satisfy the detailed balance condition
 \begin{equation}
 \label{db}
 p(x,y;\eta)e^{\sum_{\{R\mid x\in R\ {\rm or}\ y\in R\}}J_R\;\chi_R(\eta)}
 = p(y,x;\eta^{x,y})
  e^{\sum_{\{R\mid x\in R\ {\rm or}\ y\in R\}}J_R\;\chi_R(\eta^{x,y})}
 \end{equation}
 then $\nu$ is reversible for the process with rates $p$.
 \end{thm}

 \proof We want to compare the measures $\nu_N$ and $\nu_M$, where $N<M$. 
For $\eta\in Y_N$ we let $\eta'\in Y_M$ be the configuration which agrees
with $\eta$ in $T_N$ and with $\eta^H$ in $T_M\setminus T_N$, and for
$\zeta\in Y_M$ we let $\hat\zeta\in Y_N$ be the restriction of $\zeta$ to
$T_N$; thus $\hat\zeta'\equiv(\hat\zeta)'\in T_M$.  Now fix $\zeta\in Y_M$,
let $S=\{x\mid \zeta(x)\ne\hat\zeta'(x)\}$, and set $S_+=S\cap\{x\ge1\}$,
$S_-=S\cap\{x\le0\}$.  Then
 \begin{eqnarray}
  H_M(\zeta)&
   =&H_M(\hat\zeta')-2\sum_{x\in S_+}J_{\{x\}}+2\sum_{x\in S_-}J_{\{x\}}
    +\sum_{R\cap T_M\ne\emptyset\atop |R|\ge2}
      J_R[\chi_R(\zeta^*)-\chi_R(\hat\zeta'{}^*)]\nonumber\\
  &\ge&H_M(\hat\zeta')-2\sum_{x\in S_+}(J_{\{x\}}+K)
    +2\sum_{x\in S_-}(J_{\{x\}}-K).
 \end{eqnarray}
 Thus if $\eta\in Y_N$,
 \begin{equation}
  e^{-H_M(\eta')}
 \le\sum_{\hat\zeta=\eta}e^{-H_M(\zeta)}
  \le e^{-H_M(\eta')}
  \prod_{x=N+1}^M \Bigl(1+e^{2(J_{\{x\}}+K)}\Bigr)
  \prod_{x=-N}^{M-1} \Bigl(1+e^{2(-J_{\{x\}}+K)}\Bigr).
 \end{equation}
 Since the infinite products $\prod_{x\ge1}(1+e^{2(J_{\{x\}}+K)})$ and 
$\prod_{x\le0}(1+e^{2(-J_{\{x\}}+K)})$ converge by \reff{conv}, we have for
any $\epsilon>0$,
 \begin{equation}
  \label{eps}
  e^{-H_M(\eta')}
 \le\sum_{\hat\zeta=\eta}e^{-H_M(\zeta)}
  \le e^{-H_M(\eta')}(1+\epsilon),
 \end{equation}
 when $N$ is sufficiently large, uniformly in $M$. 

Now suppose that $A\subset\{0,1\}^\Z$ is such that $\one A(\eta)$ depends on
$\eta$ only through the variables $\eta(x)$ for a finite number of
sites---say for $x\in T_L$.  Since for $\eta\in Y_N$, $H_N(\eta)-H_M(\eta')$
is independent of $\eta$,
 \begin{equation}
  \nu_{N}(A)={\sum_{\eta\in Y_N,\eta^*\in A}e^{-H_{N}(\eta)}\over
     \sum_{\eta\in Y_N}e^{-H_{N}(\eta)}}
  ={\sum_{\eta\in Y_N,\eta^*\in A}e^{-H_{M}(\eta')}\over
     \sum_{\eta\in Y_N}e^{-H_{M}(\eta')}},
 \end{equation}
 and with \reff{eps} this implies that if $N\ge L$,
 \begin{equation}
  (1+\epsilon)^{-1}\nu_M(A)\le\nu_N(A)\le(1+\epsilon)\nu_M(A).
 \end{equation}
 Hence $\lim_{N\to\infty}\nu_N(A)$ exists, so that $\nu$ exists.  Similarly,
if $B\subset\{0,1\}^\Z$  is the event that  $\eta(x)=\eta^H(x)$ for
 $x\notin T_N$ then
 $\nu_M(B)=Z_M^{-1}\sum_{\eta\in Y_N}\exp(-H_M(\eta'))\ge(1+\epsilon)^{-1}$
by \reff{eps}, so that $\nu$ is a blocking measure. 
 
The measure $\nu$ is reversible for the process with rates $p$ if for any
continuous $f$ defined on $\{0,1\}^\Z$ and any $x,y\in\Z$,
 \begin{equation}
  \label{rev}
 \int p(x,y;\eta)[f(\eta^{x,y})-f(\eta)]\,d\nu = 0;
 \end{equation}
 see the proof of the analogous result for stochastic Ising models in
\cite{[Ligg]}.  But this integral may be calculated to arbitrary accuracy by
replacing $\nu$ with $\nu_N$ for suitably large $N$ (here continuity of $p$
in $\eta$ is needed), and the fact that the integral with respect to $\nu_N$
vanishes is an immediate consequence of \reff{db}.  \square

 \end{document}